\documentclass[12pt]{article}
\usepackage{amssymb}
\usepackage{latexsym,bm}
\usepackage{graphicx}
\usepackage{amsmath}
\usepackage{mathrsfs}

\date{}
\begin{document}
\title{The maximum degree of a minimally hamiltonian-connected graph\footnote{E-mail address:
\tt zhan@math.ecnu.edu.cn}}
\author{\hskip -10mm Xingzhi Zhan\\
{\hskip -10mm \small Department of Mathematics, East China Normal University, Shanghai 200241, China}}\maketitle
\begin{abstract}
 We determine the possible maximum degrees of a minimally hamiltonian-connected graph with a given order. This answers a question posed by
 Modalleliyan and Omoomi in 2016. We also pose two unsolved problems.
\end{abstract}

\vskip 8mm

We consider finite simple graphs and follow the book [5] for terminology and notations. The {\it order} of a graph is its number
of vertices, and the {\it size} is its number of edges. We denote by $V(G)$ and $E(G)$ the vertex set and edge set of a graph $G$ respectively.
For two graphs $G$ and $H,$ $G\vee H$ denotes the {\it join} of $G$ and $H,$ which is obtained from the disjoint union $G+H$
by adding edges joining every vertex of $G$ to every vertex of $H.$ $K_n$ and $C_n$ denote the complete graph of order $n$ and the cycle of order $n$
respectively. The {\it wheel} of order $n,$ denoted by $W_n,$ is the graph $K_1\vee C_{n-1}.$

A graph is called {\it hamiltonian-connected} if between any two distinct vertices there is a Hamilton path. Obviously, any hamiltonian-connected graph
of order at least $4$ is $3$-connected, and hence has minimum degree at least $3.$

{\bf Definition.} A hamiltonian-connected graph $G$ is said to be {\it minimally hamiltonian-connected} if for every edge $e\in E(G),$
the graph $G-e$ is not hamiltonian-connected.

Clearly every hamiltonian-connected graph contains a minimally hamiltonian-connected spanning subgraph. Concerning the maximum degree of a minimally hamiltonian-connected graph with a given order, Modalleliyan and Omoomi [2] proved the following results: (1) The maximum degree of any minimally hamiltonian-connected graph of order $n$ is not equal to $n-2;$ (2) the wheel $W_n$ is the only minimally hamiltonian-connected graph of order $n$ with maximum degree $n-1;$ (3) for every integer $n\ge 6$ and any integer $\Delta$ with $\lceil n/2\rceil \le\Delta\le n-3,$ there exists a minimally hamiltonian-connected graph
of order $n$ with maximum degree $\Delta.$ They [2] posed the question of whether for $\Delta$ in the range $3\le \Delta < \lceil n/2\rceil,$ there exists a minimally hamiltonian-connected graph of order $n$ with maximum degree $\Delta.$ In this note we answer the question affirmatively.

{\bf Theorem 1.}  {\it Let $n\ge 4$ be an integer. There exists a minimally hamiltonian-connected graph of order $n$ with maximum degree $\Delta$
if and only if $3\le\Delta\le n-1$ and $\Delta\ne n-2,$ where $\Delta=3$ occurs only if $n$ is even.}

{\bf Proof.} Suppose there exists a minimally hamiltonian-connected graph $G$ of order $n$ with maximum degree $\Delta.$ Then $G$ is $3$-connected, implying
that $\Delta\ge\delta\ge 3$ where $\delta$ denotes the minimum degree of $G.$  Modalleliyan and Omoomi [2] proved that $\Delta\ne n-2.$ If $\Delta=3,$ then $G$ is cubic and hence its order $n$ is even.

Conversely suppose $3\le\Delta\le n-1$ and $\Delta\ne n-2,$ and when $\Delta=3,$ $n$ is even. We will construct a minimally hamiltonian-connected graph of order $n$ with maximum degree $\Delta.$ If $\Delta=n-1,$ the wheel graph $W_n$ is a minimally hamiltonian-connected graph of order $n$ with maximum degree $n-1.$
Next suppose $3\le\Delta\le n-3.$ We will distinguish the two cases when $n-\Delta$ is odd and when $n-\Delta$ is even. For symbols such as $x_i$ below, if $i$ exceeds its valid range, then $x_i$ does not appear.

\vskip 3mm
\centerline{Case 1. $n-\Delta$ is odd}

We define a graph $G(n,\Delta)$ as follows. Denote $k=\Delta-2$ and $s=(n-\Delta+1)/2.$ We have $k\ge 1$ and $s\ge 2.$
$$
V(G(n,\Delta))=\{x_1,x_2,\ldots,x_k\}\cup\{y_1,y_2,\ldots,y_s\}\cup\{z_1,z_2,\ldots,z_{s+1}\},
$$
\begin{align*}
E(G(n,\Delta))&=\{x_ix_{i+1}|\,i=1,\ldots,k-1\}\cup\{y_iy_{i+1}|\,i=1,\ldots,s-1\}\cup\{z_iz_{i+1}|\,i=1,\ldots,s\}\\
&\cup\{y_1x_i|\,i=1,\ldots,k\}\cup\{y_iz_i|\,i=1,\ldots,s\}\cup\{x_1z_1,x_kz_{s+1},y_sz_{s+1}\}.
\end{align*}
The graph $G(16,5)$ is depicted in Figure 1.
\vskip 3mm
\par
 \centerline{\includegraphics[width=2.5in]{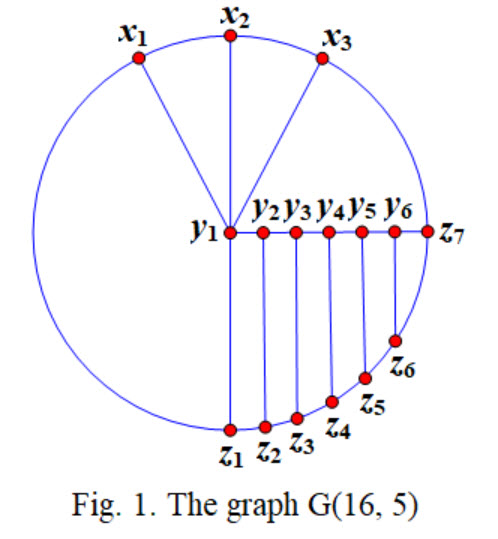}}
\par
Clearly the graph $G(n,\Delta)$ has order $n$ and maximum degree $\Delta.$
We first show that the graph $G(n,\Delta)$ is hamiltonian-connected; i.e., for any two distinct vertices $u$ and $v,$ there is a Hamilton $(u,v)$-path.
There are $8$ cases for the vertex pairs $(u,\,v)$. In each case we display a Hamilton $(u,v)$-path.

Case 1.1. $(x_i,\,x_j)$ with $1\le i<j\le k.$
\newline\indent\indent\indent\indent $x_i,x_{i+1},\ldots,x_{j-1},y_1,x_{i-1},\ldots,x_1,z_1,z_2,y_2,\ldots,z_{s+1},x_k,\ldots,x_j,$
\newline\indent\indent\indent\indent where when $i=1$ the string $x_{i-1},\ldots,x_1$ does not appear.

Case 1.2. $(x_i,\,y_j)$ with $1\le i\le k$ and $1\le j\le s.$
\newline\indent\indent\indent\indent $x_i,x_{i+1},\ldots,x_k,z_{s+1},\ldots,y_{j+1},z_{j+1},z_j,\ldots,z_1,x_1,\ldots,x_{i-1},y_1,\ldots,y_j,$
\newline\indent\indent\indent\indent where when $i=1$ the string $x_1,\ldots,x_{i-1}$ does not appear.

Case 1.3. $(x_i,\,z_j)$ with $1\le i\le k$ and $1\le j\le s.$
\newline\indent\indent\indent\indent $x_i,x_{i+1}\ldots,x_k,z_{s+1},\ldots,z_{j+1},y_{j+1},y_j,\ldots,y_1,x_{i-1},\ldots,x_1,z_1,\ldots,z_j,$
\newline\indent\indent\indent\indent where when $i=1$ the string $x_{i-1},\ldots,x_1$ does not appear.

Case 1.4. $(x_i,\,z_{s+1})$ with $1\le i\le k.$
\newline\indent\indent\indent\indent $x_i,x_{i+1},\ldots,x_k,y_1,x_{i-1},\ldots,x_1,z_1,z_2,\ldots,z_{s+1},$
\newline\indent\indent\indent\indent where when $i=1$ the string $x_{i-1},\ldots,x_1$ does not appear.

Case 1.5. $(y_i,\,y_j)$ with $1\le i<j\le s.$
\newline\indent\indent\indent\indent $y_i,y_{i+1},\ldots,y_{j-1},z_{j-1},\ldots,z_i,z_{i-1},\ldots,x_1,x_2,\ldots,x_k,z_{s+1},\ldots,z_j,y_j.$

Case 1.6. $(y_i,\,z_j)$ with $1\le i<j\le s+1.$
\newline\indent\indent\indent\indent $y_i,y_{i+1},\ldots,y_{j-1},z_{j-1},\ldots,z_i,z_{i-1},\ldots,x_1,x_2,\ldots,x_k,z_{s+1},\ldots,y_j,z_j.$

Case 1.7. $(y_i,\,z_j)$ with $1\le j\le i\le s.$
\newline\indent\indent\indent\indent $y_i,y_{i-1},\ldots,y_j,y_{j-1},z_{j-1},\ldots,x_1,x_2,\ldots,x_k,z_{s+1},\ldots,y_{i+1},z_{i+1},z_i,\ldots,z_j.$

Case 1.8. $(z_i,\,z_j)$ with $1\le i<j\le s+1.$
\newline\indent\indent\indent\indent $z_i,z_{i+1},\ldots,z_{j-1},y_{j-1},\ldots,y_i,y_{i-1},z_{i-1},\ldots,x_1,x_2,\ldots,x_k,z_{s+1},\ldots,z_j.$

We have shown that $G(n,\Delta)$ is hamiltonian-connected. Recall that any hamiltonian-connected graph of order at least $4$ is $3$-connected and hence
has minimum degree at least $3.$ Note that every edge of $G(n,\Delta)$ has one endpoint of degree $3.$ Thus for any $e\in E(G(n,\Delta)),$
$G(n,\Delta)-e$ has a vertex of degree $2,$ implying that it is not hamiltonian-connected. This completes the proof that $G(n,\Delta)$ is minimally hamiltonian-connected.

\vskip 3mm
\centerline{Case 2. $n-\Delta$ is even}

We define a graph $H(n,\Delta)$ as follows. Denote $k=\Delta-1$ and $s=(n-\Delta-2)/2.$ Since $n-\Delta$ is even and $\Delta\le n-3,$ we have $k\ge 3$
and $s\ge 1.$
$$
V(H(n,\Delta))=\{x\}\cup \{y_1,y_2,\ldots,y_k\}\cup\{z_0,z_1,z_2,\ldots,z_s\}\cup\{w_1,w_2,\ldots,w_{s+1}\},
$$
\begin{align*}
E(H(n,\Delta))&=\{y_iy_{i+1}|\,i=1,\ldots,k-1\}\cup \{z_iz_{i+1}|\,i=0,1,\ldots,s-1\}\cup\{w_iw_{i+1}|\,i=1,\ldots,s\}\\
              &\cup\{xy_i|\,i=1,\ldots,k\}\cup\{z_iw_i|\,i=1,\ldots,s\}\cup\{xz_1,y_1z_0,z_0w_1,y_kw_{s+1},z_sw_{s+1}\}.
\end{align*}
The graph $H(17,5)$ is depicted in Figure 2.
\vskip 3mm
\par
 \centerline{\includegraphics[width=2.5in]{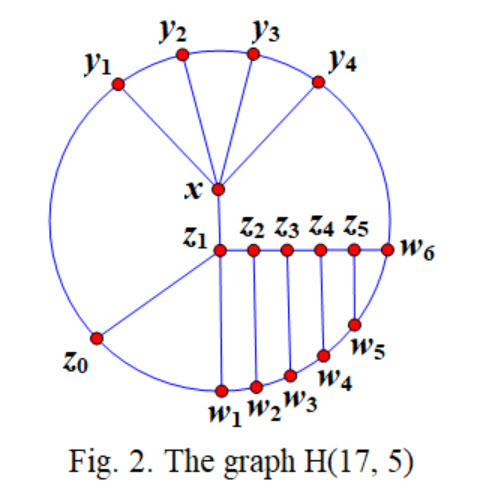}}
\par
Clearly the graph $H(n,\Delta)$ has order $n$ and maximum degree $\Delta.$
We first show that the graph $H(n,\Delta)$ is hamiltonian-connected; i.e., for any two distinct vertices $u$ and $v,$ there is a Hamilton $(u,v)$-path.
There are $16$ cases for the vertex pairs $(u,\,v)$. In each case we display a Hamilton $(u,v)$-path.

Case 2.1. $(y_1,\,y_j)$ with $2\le j\le k.$
\newline\indent\indent\indent\indent $y_1,\ldots,y_{j-1},x,z_1,z_0,w_1,\ldots,w_{s+1},y_k,\ldots,y_j.$

Case 2.2. $(y_i,\,y_j)$ with $2\le i<j\le k.$
\newline\indent\indent\indent\indent $y_i,y_{i+1},\ldots,y_{j-1},x,y_{i-1},\ldots,y_1,z_0,\ldots,w_{s+1},y_k,\ldots,y_j.$

Case 2.3. $(y_1,\,x).$
\newline\indent\indent\indent\indent $y_1,z_0,\ldots,w_{s+1},y_k,\ldots,y_2,x.$

Case 2.4. $(y_i,\,x)$ with $2\le i\le k.$
\newline\indent\indent\indent\indent $y_i,y_{i+1},\ldots,y_k,w_{s+1},\ldots,z_0,y_1,\ldots,y_{i-1},x.$

Case 2.5. $(y_i,\,z_j)$ with $1\le i\le k-1$ and $0\le j\le s.$
\newline\indent\indent\indent\indent $y_i,y_{i-1},\ldots,y_1,x,y_{i+1},\ldots,y_k,w_{s+1},\ldots,z_{j+1},w_{j+1},w_j,\ldots,w_1,z_0,\ldots,z_j.$

Case 2.6. $(y_k,\,z_0).$
\newline\indent\indent\indent\indent $y_k,y_{k-1},\ldots,y_1,x,z_1,\ldots,z_s,w_{s+1},\ldots,w_1,z_0.$

Case 2.7. $(y_k,\,z_j)$ with $1\le j\le s.$
\newline\indent\indent\indent\indent $y_k,y_{k-1},\ldots,y_2,x,y_1,z_0,\ldots,z_{j-1},w_{j-1},w_j,\ldots,w_{s+1},z_s,\ldots,z_j.$

Case 2.8. $(y_1,\,w_j)$ with $1\le j\le s+1.$
\newline\indent\indent\indent\indent $y_1,z_0,w_1,\ldots,w_{j-1},z_{j-1},\ldots,z_1,x,y_2,\ldots,y_k,w_{s+1},\ldots,w_j,$
\newline\indent\indent\indent\indent where when $j=1$ the string $w_1,\ldots,w_{j-1},z_{j-1},\ldots,z_1$ means $z_1.$

Case 2.9. $(y_i,\,w_j)$ with $2\le i\le k$ and $1\le j\le s+1.$
\newline\indent\indent\indent\indent $y_i,y_{i+1},\ldots,y_k,x,y_{i-1},\ldots,y_1,z_0,\ldots,w_{j-1},z_{j-1},z_j,\ldots,z_s,w_{s+1},\ldots,w_j.$

Case 2.10. $(x,\,z_j)$ with $0\le j\le s.$
\newline\indent\indent\indent\indent $x,y_1,\ldots,y_k,w_{s+1},\ldots,z_{j+1},w_{j+1},w_j,\ldots,w_1,z_0,\ldots,z_j.$

Case 2.11. $(x,\,w_j)$ with $1\le j\le s+1.$
\newline\indent\indent\indent\indent $x,y_k,\ldots,y_1,z_0,\ldots,w_{j-1},z_{j-1},z_j,\ldots,z_s,w_{s+1},\ldots,w_j.$

Case 2.12. $(z_i,\,z_j)$ with $0\le i<j\le s.$
\newline\indent\indent\indent\indent $z_i,w_i,w_{i-1},\ldots,z_0,y_1,x,y_2,\ldots,y_k,w_{s+1},\ldots,z_{j+1},w_{j+1},w_j,\ldots,w_{i+1},z_{i+1},\ldots,z_j.$

Case 2.13. $(z_0,\,w_j)$ with $1\le j\le s+1.$
\newline\indent\indent\indent\indent $z_0,w_1,\ldots,w_{j-1},z_{j-1},\ldots,z_1,x,y_1,\ldots,y_k,w_{s+1},\ldots,w_j,$
\newline\indent\indent\indent\indent where when $j=1$ the string $w_1,\ldots,w_{j-1},z_{j-1},\ldots,z_1$ means $z_1.$

Case 2.14. $(z_i,\,w_j)$ with $1\le i<j\le s+1.$
\newline\indent\indent\indent  $z_i,z_{i+1},\ldots,z_{j-1},w_{j-1},w_{j-2},\ldots,w_i,w_{i-1},\ldots,z_0,y_1,x,y_2,\ldots,y_k,w_{s+1},\ldots,z_j,w_j.$

Case 2.15. $(z_i,\,w_j)$ with $1\le j\le i\le s.$
\newline\indent\indent\indent  $z_i,z_{i-1},\ldots,z_j,z_{j-1},w_{j-1},\ldots,z_0,y_1,x,y_2,\ldots,y_k,w_{s+1},\ldots,z_{i+1},w_{i+1},\ldots,w_j.$

Case 2.16. $(w_i,\,w_j)$ with $1\le i<j\le s+1.$
\newline\indent\indent\indent $w_i,w_{i+1},\ldots,w_{j-1},z_{j-1},\ldots,z_i,z_{i-1},w_{i-1},\ldots,z_0,y_1,x,y_2,\ldots,y_k,w_{s+1},\ldots,z_j,w_j.$

Thus we have shown that $H(n,\Delta)$ is hamiltonian-connected. Recall that any hamiltonian-connected graph of order at least $4$ is $3$-connected and hence
has minimum degree at least $3.$ Since the graph $H(n,\Delta)-xz_1$ has connectivity $2,$ it is not hamiltonian-connected. For every edge
$e\in E(H(n,\Delta))$ with $e\ne xz_1,$ $e$ has one endpoint of degree $3.$ Therefore $H(n,\Delta)-e$ has a vertex of degree $2,$ implying that it is not
hamiltonian-connected. This completes the proof that $H(n,\Delta)$ is minimally hamiltonian-connected, and the theorem is proved. $\Box$

{\bf Remark.} The graphs $G(n,\Delta)$ and $H(n,\Delta)$ constructed in the above proof of Theorem 1 have degree sequences $\Delta,3,3,\ldots,3$ and $\Delta,4,3,\ldots,3$ respectively, and hence they have the minimum possible sizes among all graphs of order $n$ with maximum degree $\Delta$
and minimum degree at least $3$ in the two cases when $n-\Delta$ is odd and when $n-\Delta$ is even respectively. The graph constructed in [2] for $\Delta$ in the range $\lceil n/2\rceil\le\Delta\le n-3$ has degree sequence $\Delta, n-\Delta, 3,\ldots,3.$

Finally we pose two unsolved problems.

{\bf Problem 1.} Let $n\ge 4$ be a given integer. What are the possible values of the minimum degree of a minimally hamiltonian-connected graph of order $n$?

A computer search shows that every minimally hamiltonian-connected graph of order $n$ with $4\le n\le 10$ has minimum degree $3.$ The author does not know of
an example of a minimally hamiltonian-connected graph with minimum degree at least $4.$ The following easier problem is of a more basic nature.

{\bf Problem 2.} Does there exist a minimally hamiltonian-connected graph with minimum degree at least $4$?

There are some sufficient conditions for hamiltonian-connected graphs; for recent ones see [1], [3] and [4]. But very little is known about necessary conditions.
Restrictions on the maximum or minimum degree of a minimally hamiltonian-connected graph may be viewed as necessary conditions for this smaller class of
graphs.

\vskip 5mm
{\bf Acknowledgement.} This research  was supported by the NSFC grant 11671148 and Science and Technology Commission of Shanghai Municipality (STCSM) grant 18dz2271000.

\end{document}